\documentclass[12pt]{article}
\newtheorem{theo}{Theorem}
\newtheorem{prop}{Proposition}

\begin{document}
\newpage
\pagestyle{empty}
\setcounter{page}{0}
%
\newfont{\twelvemsb}{msbm10 scaled\magstep1}
\newfont{\eightmsb}{msbm8} \newfont{\sixmsb}{msbm6} \newfam\msbfam
\textfont\msbfam=\twelvemsb \scriptfont\msbfam=\eightmsb
\scriptscriptfont\msbfam=\sixmsb \catcode`\@=11
\def\Bbb{\ifmmode\let\next\Bbb@\else \def\next{\errmessage{Use
      \string\Bbb\space only in math mode}}\fi\next}
\def\Bbb@#1{{\Bbb@@{#1}}} \def\Bbb@@#1{\fam\msbfam#1}
\newfont{\twelvegoth}{eufm10 scaled\magstep1}
\newfont{\tengoth}{eufm10} \newfont{\eightgoth}{eufm8}
\newfont{\sixgoth}{eufm6} \newfam\gothfam
\textfont\gothfam=\twelvegoth \scriptfont\gothfam=\eightgoth
\scriptscriptfont\gothfam=\sixgoth \def\frak{\frak@}
\def\frak@#1{{\fam\gothfam{{#1}}}} \def\frak@@#1{\fam\gothfam#1}
\catcode`@=12
%
%
%
\def\CC{{\Bbb C}}
\def\NN{{\Bbb N}}
\def\QQ{{\Bbb Q}}
\def\RR{{\Bbb R}}
\def\ZZ{{\Bbb Z}}
\def\cA{{\cal A}}          \def\cB{{\cal B}}          \def\cC{{\cal C}}
\def\cD{{\cal D}}          \def\cE{{\cal E}}          \def\cF{{\cal F}}
\def\cG{{\cal G}}          \def\cH{{\cal H}}          \def\cI{{\cal I}}
\def\cJ{{\cal J}}          \def\cK{{\cal K}}          \def\cL{{\cal L}} 
\def\cM{{\cal M}}          \def\cN{{\cal N}}          \def\cO{{\cal O}}
\def\cP{{\cal P}}          \def\cQ{{\cal Q}}          \def\cR{{\cal R}} 
\def\cS{{\cal S}}          \def\cT{{\cal T}}          \def\cU{{\cal U}}
\def\cV{{\cal V}}          \def\cW{{\cal W}}          \def\cX{{\cal X}}
\def\cY{{\cal Y}}          \def\cZ{{\cal Z}}
\def\qed{\hfill \rule{5pt}{5pt}}
\def\th{\mbox{\footnotesize th}}
\newenvironment{result}{\vspace{.2cm} \em}{\vspace{.2cm}}
%
%
\newcommand{\norm}[1]{{\protect\normalsize{#1}}}
\newcommand{\LAP}
{{\small E}\norm{N}{\large S}{\Large L}{\large A}\norm{P}{\small P}}
\newcommand{\sLAP}{{\scriptsize E}{\footnotesize{N}}{\small S}{\norm L}$
${\small A}{\footnotesize{P}}{\scriptsize P}}
\def\logolapin{
  \raisebox{-1.2cm}{\epsfbox{/lapphp8/keklapp/ragoucy/paper/enslapp.ps}}}
\def\logolight{{\bf{{\large E}{\Large N}{\LARGE S}{\huge L}{\LARGE
        A}{\Large P}{\large P}} }}
\def\logoenslapp{\logolight}
%
%
%
\hbox to \hsize{
\hss
\begin{minipage}{5.2cm}
  \begin{center}
    {\bf Groupe d'Annecy\\ \ \\
      Laboratoire d'Annecy-le-Vieux de Physique des Particules}
  \end{center}
\end{minipage}
\hfill
\logoenslapp
\hfill
\begin{minipage}{4.2cm}
  \begin{center}
    {\bf Groupe de Lyon\\ \ \\
      {\'E}cole Normale Sup{\'e}rieure de Lyon}
  \end{center}
\end{minipage}
\hss}

\vspace {.3cm}
\centerline{\rule{12cm}{.42mm}}
%
\vfill
\vfill
\begin{center}

  {\LARGE {\bf {\sf Scasimir operator, Scentre 
        and \\[2mm]
        Representations of $\cU_q(osp(1|2))$ }}}
 \\[1cm]

\vfill

{\large D. Arnaudon$^{\ddagger,}$
\footnote{\ arnaudon@lapp.in2p3.fr
  
\indent
  \ $^2$ bauer@spht.saclay.cea.fr

  \indent
  \ $^3$ Permanent address: Service de Physique Th{\'e}orique,
  C.E.A. Saclay, F-91191, Gif-sur-Yvette, France. 

  \indent
  \ $^4$ URA 14-36 du CNRS, associ{\'e}e {\`a} l'E.N.S. de Lyon et {\`a}
  l'Universit{\'e} de Savoie.}
and M. Bauer$^{\ddagger,2,3}$}

\vfill

{\large $^\ddagger$}{\em \LAP$^4$, Chemin de Bellevue BP 110,
74941 Annecy-le-Vieux Cedex, France.}

\end{center}

\vfill

\begin{abstract}

A bosonic operator of $\cU_q(osp(1|2))$ that anticommutes with the
fermionic generators appears to be useful to describe the relations in
the centre of $\cU_q(osp(1|2))$ for $q$ a root of unity (in the
unrestricted specialisation). 
As in the classical case, it also
simplifies the classification of finite dimensional irreducible
representations. 

\end{abstract}

\vfill
\vfill

\rightline{\LAP --A--591/96}
\rightline{SPhT--96--xxx}
\rightline{q-alg/9605020}
\rightline{May 96}

\newpage
\pagestyle{plain}
\setcounter{footnote}{0}

\section{Motivations\label{sect:motivations}}

The centre of the universal enveloping algebra of a semi simple Lie
algebra is a free polynomial algebra that can be described quite
explicitly with the Harish-Chandra homomorphism. This construction
can be extended to the quantum deformation of a semi simple Lie
algebra when the deformation parameter $q$ is not 
a root of unity \cite{RosHC}, leading to deformed Casimir operators.

But when $q$ is a root of unity\footnote{We work in the unrestricted
specialisation, so we do not introduce divided powers of the
generators.}, say $q^l=1$,
the Casimir operators do 
not exhaust the centre anymore. For example if $l$ is odd, the
$l^{\th}$ powers of the generators are central. Together with the
Casimir operators, they satisfy polynomial relations
\cite{Lusztigmod,DK,KerDipl,ABrelations,AACuqsln}. 
Hence from the
point of view of algebraic geometry, the centre is a non trivial
affine variety. From a physical point of view, it is desirable to know
the equations of this variety. This is because to exploit efficiently
the symmetries of a model one needs a complete set of commuting
observables to label the states. If the symmetry of the model is a
quantum group associated to a semi simple Lie algebra, the centre
together with a Cartan subalgebra gives
 such a complete set. 

For semi simple superalgebras, there is no strict analogue of the
Harish-Chandra construction 
in general \cite{KacSuperRep}. 
The centre of the universal enveloping algebra does not
have to be a free polynomial algebra. But the Casimir operators can
still be extended to the quantum deformation and exhaust the centre
when $q$ is not a root of unity. 

When $q$ is a root of unity, the discussion of the Lie algebra case applies. 

The relations in the centre of a quantum deformation of a semi simple
(Lie or super) algebra have been computed only for particular
examples \cite{KerDipl,ABrelations,AACuqsln,AAB}. 
All the known cases seem to share some nice
features that we 
briefly recall in section  \ref{sect:observations}. 

The study of $\cU_q(osp(1|2))$ that we present below has led us to
reconsider those general features. In fact, trying to conciliate
them with the case at hand, we were led to 
use the following simple structure: the Casimir
operator of $\cU_q(osp(1|2))$ \cite{KulRes}  is a perfect square, and
the square 
root, despite its bosonic character, anti-commutes with fermions and
commutes with bosons, so we decided to call it a Scasimir
operator. 
This operator, first written in \cite{Lesniewski} is the
$q$-deformation of a classical operator 
introduced in \cite{PaisRitt, Pinc}. 
The previously observed general features of the relations in
the centre satisfied by the Casimir operators now involve the Scasimir
operator. The possible generalisation to  $\cU_q(osp(1|2n))$ involves
a detailed study of the Scasimir operator in the non deformed
case. This study is presented elsewhere \cite{ABFscasimirs}. The
importance of the 
Scasimir operator to classify irreducible finite dimensional
irreducible representations of $\cU_q(osp(1|2))$ is also emphasised. 

We note that the existence of periodic irreducible representations
(for $q$ a root of unity, in the unrestricted case)
imply the existence of primitive ideals that are not annihilators of
irreducible quotients of Verma modules (all of them being annihilated
by a common finite power of all the raising generators). This differs
from the case of classical (non quantum) (super)algebras
\cite{Duflo,MussonPrim}. 

\section{Notations\label{sect:notations}}

In this paper, $q$ is a complex number such that $q^2 \neq 0,1$ and
$l$ is an integer larger than $2$. We also need an indeterminate $t$
to build generating functions.

We denote by $q^{1/2}$ a fixed square root of $q$ and by $q'$ the
opposite of $q$. We set $\eta = (q^{1/2}+q^{-1/2}) (q-q^{-1})$.

In the following, we shall often be interested in the case when $q$ is
a primitive $l^{\th}$-root of unity. We use $l'$ to denote the order of
$q'$ and $L$ to denote the smallest even multiple of $l$ (that is, $l$ if
$l$ is even and $2l$ if $l$ is odd). The integer $l'$ is $l/2$ if $l$
is twice an odd integer and $L$ otherwise. The map $l \rightarrow l'$
is one to one. 

\section{Definitions\label{sect:definitions}}

The algebra $\cU_q(osp(1|2))$ is the unital associative algebra with 
generators $e,f,k,k^{-1}$ and relations
\begin{eqnarray}
  && k e k^{-1} = q e \;,   \qquad\qquad
  k f k^{-1} = q^{-1} f\;,   \\
  && e f + f e = \frac{k-k^{-1}}{q-q^{-1}}\;, \qquad\qquad k k^{-1} =
  k^{-1} k =1\;.
\end{eqnarray}

The $\ZZ _2$ grading $d(e)=d(f)=1$, $d(k)=0$ is compatible with the
relations, hence has a unique extension to a grading of
$\cU_q(osp(1|2))$. We write $\cU_q(osp(1|2))=\cU_q(osp(1|2))_0 \oplus 
\cU_q(osp(1|2))_1$ and refer to elements in the first (resp. second)
summand as bosons (resp. fermions).  As usual $\cU_q(osp(1|2))_0$ is a
subalgebra of $\cU_q(osp(1|2))$. 

The algebra $\cU_q(osp(1|2))$ has a number of obvious auto-morphisms. For
 instance the scale change $(e,f,k,k^{-1}) \rightarrow (ae,a^{-1}f,k,k^{-1})$
(where $a$ is a nonzero complex number), the signed permutation  
$(e,f,k,k^{-1}) \rightarrow (-f,e,k^{-1},k)$ and combinations thereof
induce auto-morphisms.

Of central interest for us is $S$, the Scasimir operator. It is the
boson defined by
\begin{equation} \label{eq:scas}
  S = q^{1/2} k - q^{-1/2} k^{-1} - \eta fe
\end{equation}
where $\eta = (q^{1/2}+q^{-1/2}) (q-q^{-1})$.

One can check directly the remarkable fact that
\begin{prop} [Le\'sniewski \cite{Lesniewski}] 
  The operator $S$ anti-commutes with
  fer\-mions and commutes with bosons.
  Its square $S^2$ is nothing but the standard Casimir
  operator $C$ up to a constant. 
\end{prop}


Explicitly, 
\begin{equation}
  C = qk^2 + q^{-1} k^{-2} +(q-q^{-1})^2 (qk+q^{-1}k^{-1})fe
      - (q-q^{-1})^2 (q+2+q^{-1}) f^2 e^2 
\end{equation}
is related to $S$ by 
\begin{equation}
  S^2 + 2 = C.
\end{equation}

The existence of $S$ is not a byproduct of the quantum deformation. Up
to an overall factor, $S$ has a well defined classical limit when $q 
\rightarrow 1$. Set $k=q^h$ and get
\begin{equation}
 S_{class}=h-2fe+\frac{1}{2}=ef-fe+\frac{1}{2}.
\end{equation}
A normalized version of this operator has first been introduced in
\cite{PaisRitt} where it was used as grading operator on
representations. 
This operator has also been defined in \cite{Pinc}, where it is proved
that it generates a minimal primitive ideal without Lie algebra
analogue since it is not generated by its intersection with the
centre. 
The operator $S$ is also found as part of the image of the
quadratic Casimir operator of $\cU_q(sl(2))$ by the application
defined in \cite{Zhangd}.

\section{Commutations \label{sect:commutations}}

Recall that $q'$ is the opposite of $q$. Let $m$ be a positive integer. We
define $\varepsilon(m)$ to be $1$ if $m=0,1 \; \bmod 4$ and $-1$ otherwise.

In the sequel, we shall need explicit (anti)commutations relations
between powers of $e$ and $f$. A recursion argument shows that 
\begin{equation} \label{eq:comemf}
(q-q^{-1})(f^me+(-1)^{m-1}ef^m)=f^{m-1}\left(\frac{q'^{-m}-1}{q'^{-1}-1}k
-\frac{q'^m-1}{q'-1}k^{-1}\right).
\end{equation}
The corresponding equation when the roles of $e$ and $f$ are interchanged
can be obtained without computations using the auto\-morphism
$(e,f,k,k^{-1}) \rightarrow (-f,e,k^{-1},k)$.

We adapt a trick used by Kerler 
\cite{KerDipl} for $\cU_q(sl(2))$ to get an identity relating 
the Scasimir operator and $k$ to powers of $e$ and $f$. The formula
\begin{equation} \label{eq:scasm}
\prod_{n=0}^{m-1}(S-q'^nq^{1/2}k+q'^{-n}q^{-1/2}k^{-1})=\varepsilon(m)
(-\eta)^m f^m e^m
\end{equation}
is proved by recursion starting from (\ref{eq:scas}). To go from $m$
to $m+1$ one multiplies both sides by $f$ on the left, by $e$ on the
right and by $-\eta\; \varepsilon(m)\varepsilon(m+1)$. The right hand
side is what is expected. On the left hand side, $f$ goes through the
factors, multiplying $k$ by $-q'$, $k^{-1}$ by $-q'^{-1}$ and $S$ by
$-1$. Finally it reaches $e$ when (\ref{eq:scas}) is used to eliminate
$-\eta fe$. The signs disappear because $(-1)^m \varepsilon(m)
\varepsilon (m+1)=1$. 

\section{Foundations \label{sect:foundations}}

In the sequel, we shall use repeatedly the following

\begin{prop} The family $\{f^ae^bk^c; \; a,b\in\NN,\; c\in\ZZ\}$ is a
linear basis of $\cU_q(osp(1|2))$.
\end{prop}
 
Applying the (anti-)commutation relations, any element of
$\cU_q(osp(1|2))$ can be written as a linear combination of elements
in $\{f^ae^bk^c; \; a,b\in\NN,\; c\in\ZZ\}$. That this decomposition
is unique is the content of a (super, quantum) analogue of the
Poincar{\'e}--Birkoff--Witt theorem. To prove it in this special case,
one builds $\cU_q(osp(1|2))$ by two successive Ore extensions of the
algebra $A_0$ of Laurent polynomials in $k$. Formula (\ref{eq:comemf})
plays a crucial role. We do not give the
details because they follow closely the proof given in 
\cite{Kassel}, Chap VI, for $\cU_q(sl(2))$. Chapter I of \cite{Kassel}
contains a good introduction to Ore extensions.   

\medskip

We will now prove

\begin{prop}
  If $q$ is not a root of unity, the centre of $\cU_q(osp(1|2))$ is $\CC
  [C]$ (where $C=S^2+2$ is the Casimir operator).
\end{prop}
{\em Proof: }
Consider the commutative subalgebra $\cA$ of $\cU_q(osp(1|2))$ with 
generators $S,k,k^{-1}$. First, $\cA$ is isomorphic to $\CC
[k,k^{-1},S]$. This amounts to say that there is no polynomial
relation between $k$ and $S$. This is because nothing could compensate
the monomial of highest degree in $S$ (say $m$), the only one which
produces a 
term $f^me^m$ when written in the Poincar{\'e}--Birkoff--Witt basis. 

Applying formula (\ref{eq:scasm}) with $m=\min (a,b)$ to the monomial 
$f^ae^bk^c$ shows that $\cU_q(osp(1|2))$ is a free $\cA$-module with
basis $1,f,e,f^2,e^2,\cdots$.
In this basis, it is easy to look for central elements. Until the next
section, we assume that $q$ is not a root of unity. The adjoint
action of $k$ is diagonal, with eigenvalue
$1,q^{-1},q,q^{-2},q^2,\cdots$ for $1,f,e,f^2,e^2,\cdots$. Hence the
commutant of $k$ in $\cU_q(osp(1|2))$ is $\cA$. Now we have to look for
the commutant of $f$ and $e$ in $\cA$. For $P(k,S) \in \cA$, it is
easily shown that $P(k,S)e=eP(qk,-S)$. The commutation condition 
$P(kq,-S)=P(k,S)$ implies that $P$ has to be an even function of $S$
independent of $k$. To see this, expand $P(k,S)$ in monomials and
identify term by term : a monomial $k^iS^j$ can appear with non zero
coefficient only if $q^i(-1)^j=1$. As $q$ is not a root of unity, $i$
has to vanish and $j$ has to be even. The commutation with $f$ gives
the same condition. 

\qed

In the next section, we shall address the question of the structure of
the centre when $q$ is a root of unity. 

\section{Restrictions \label{sect:restrictions}}  
 
For the rest of the paper, $l$ is an integer larger than $2$ and $q$
is a primitive $l^{\th}$ root of unity. The integer $L$ is the smallest
even multiple of $l$ (that is, $l$ if $l$ is even and $2l$ if $l$ is
odd). The integer $l'$ is $L/2$ if $l$ is twice an odd integer and $L$
otherwise. Then $q'=-q$ is a primitive $l'^{\th}$ root of unity, and $q^L=1$. 
We shall now give a complete description of the centre in this case.

Evaluation of (\ref{eq:comemf}) for $m=L$ gives
\begin{equation} f^{L}e-ef^{L}=0.
\end{equation}
Due to auto-morphisms, the corresponding relation with $e$ and $f$
interchanged is true as well. But $F\equiv f^L$ and $E \equiv e^L$
commute with $k$. So they are central. We observe again that we do not
introduce divided powers, so $F$ and $E$ can be non-zero in representations.

We can now compute the structure of the centre.

\begin{theo}
  The centre
  of $\cU_q(osp(1|2))$ contains $\CC [k^l,k^{-l},C,F] + \CC
  [k^l,k^{-l},C,E]$ (this sum is not direct : $\CC [k^l,k^{-l},C,F]
  \cap \CC [k^l,k^{-l},C,E] = \CC [k^l,k^{-l},C]$).  
  \begin{itemize}
  \item 
    If $l$ is odd, this
    is the entire centre. 
  \item
    If $l$ is even, the cent\-re is a free module over 
    $\CC [k^l,k^{-l},C,F] + \CC [k^l,k^{-l},C,E]$ with basis $1,k^{l/2}S$.
  \end{itemize}
\end{theo}
{\em Proof:}
The commutative (sub-)algebras (of $\cU_q(osp(1|2))$) $\cA [F]$ and
$\cA [E]$ are free polynomial algebras (the argument given for $\cA$ 
in the previous section still applies). We can use $F$ and $E$ to
refine the previously obtained decomposition of $\cU_q(osp(1|2))$ as 
an $\cA$-module : $\cU_q(osp(1|2))=\cU _+ +\cU _-$ where $\cU _+$ is 
the free $\cA [F]$-module with basis $1,f,\cdots,f^{L-1}$ and $\cU _-$
the free $\cA [E]$-module with basis $1,e,\cdots,e^{L-1}$. The sum $\cU _+ 
+\cU _-$ is not direct, but $\cU _+ \cap \cU _-=\cA$. 

The monomial $e^m$ commutes with $k$ only if $m$ is a multiple of $l$
and with $S$ only if $m$ is even. The same is true for powers of $f$. 
Hence the commutant of $k$ and $S$ in $\cU_q(osp(1|2))=\cU _+ +\cU _-$
is $\cA [F] + \cA [E]$. Again, this sum is not direct and $\cA [F]
\cap \cA [E] = \cA$. We adapt the argument of the previous section. Any 
element of $\cA [F] + \cA [E]$ can be expanded in powers of $k$ and
$S$. In the expansion of a central element, a monomial $k^iS^j$ can
appear with non zero coefficient only if $q^i(-1)^j=1$. This happens
for instance if $i$ is a multiple of $l$ and $j$ is even. If $l$ is
odd, this is the only possibility. If $l$ is even, there is another
solution, namely, $i$ an odd multiple of $l/2$ and $j$ odd. 

\qed

This gives a simple explicit description  of the centre, but there are
some drawbacks related to the multiplicative structure. For instance,
although $E$ and $F$ are described as such, we still need an expression
for $FE$. This is the purpose of section \ref{sect:relations}. 

Moreover, evaluation of (\ref{eq:comemf}) for $m=l'$ gives
\begin{equation} f^{l'}e+(-1)^{l'-1}ef^{l'}=0.
\end{equation}
Due to auto-morphisms, the corresponding relation with $e$ and $f$
interchanged is true as well. Assume $l$ is twice an odd integer, so that 
$l=L=2l'$ and $l'$ is odd. Then $e^{l'}$ anti-commutes with $k$,
$k^{-1}$ and $S$. The above equation shows that it also anti-commutes 
with $f$. Analogously, $f^{l'}$ anti-commutes with $k$, $k^{-1}$,$S$ 
and $e$. So there are unexpected central elements, namely 
$e^{L/2}f^{L/2}k^{\pm L/2}$. They differ by a factor $k^l$ so only one
of those needs to be considered. It will be expressed in the standard
description in section \ref{sect:relations} as well. 

\section{Observations\label{sect:observations}}

Let us recall some standard facts. The quantum deformation $A_q$ of a
semi simple (Lie or super) algebra depends on a complex parameter
$q$. It is generated by raising operators, lowering operators, and
Cartan generators. By a choice of Poincar{\'e}--Birkoff--Witt basis (a
standard choice is to put lowering operators on the  left, Cartan
generators in the middle and raising operators on the right), we can
identify the vector spaces underlying the algebras $A_q$ for different
values of $q$, and it makes sense to talk about a family of elements
in $A_q$ depending smoothly on $q$. A Casimir operator is such a
family of central elements. When $q$ is not a root of unity, the
Casimirs span the centre of $A_q$. When $q$ is a root of unity, there
is an integer $L$, simply related to the order of $q$, such that the
$L^{\th}$ power of any generator is central. Those new central elements
do not satisfy any relation among themselves because such relations
would imply relations in the Poincar{\'e}--Birkoff--Witt basis (this is
automatic for the Lie case, for the super case, one has to be more
cautious). At the same time, the Casimirs are not independent from
those new generators : there are relations in the centre. 

In cases treated up to now (i.e. $\cU_q(sl(2))$ in \cite{KerDipl},
$\cU_q(sl(N))$ in \cite{ABrelations,AACuqsln}, 
$\cU_q(sl(2|1))$ in \cite{AAB}), the
following general features have been observed.  
After a choice of Poincar{\'e}--Birkoff--Witt basis, elements of $A_q$
are identified with (non commuting) polynomials, and the notion of 
substitution is well defined. For instance, if $C$ is a Casimir
element, one can replace all the coefficients and generators by their
$L^{\th}$ power and get a new (central because expressed in terms of
central powers of the generators) element $C^{(L)}$. It turns out that
$C^{(L)}$ is a polynomial in the Casimirs. This is already
remarkable. Moreover, this polynomial has a simple description. It is
a generalised Chebychev polynomial : substitution of $0$ for the
raising and lowering operators gives a restriction involving only Cartan
generators, which is enough to compute the desired relation.

This general setting is a little bit abstract. Our point is that it does not
work as it stands for $\cU_q(osp(1|2))$. However only a minor
modification is needed. It is the failure of the general philosophy
that motivated us to look for the Scasimir. Once we use the Scasimir
instead of the Casimir, the above construction works. So a glance at
the section \ref{sect:relations} will provide a concrete 
example of the general facts that we just outlined.

\section{Functions\label{sect:functions}}

Our subsequent study makes essential use of families of polynomials
which we describe now. 

If $u$ is an indeterminate and $m$ a positive integer, we claim that 
$u^m+(-u^{-1})^m$ is a monic polynomial of order $n$ in $S=u-u^{-1}$
of parity $(-1)^m$, in fact a Chebychev polynomial. This comes from
trigonometric identities, but we prefer to consider the generating function
\begin{equation}
\sum_{m \geq 0} t^m
(u^m+(-u^{-1})^m)=\frac{1}{1-tu}+\frac{1}{1+tu^{-1}}= 
\frac{2-tS}{1-tS-t^2}=\sum_{m \geq 0} t^m P_m(S).
\end{equation}
The polynomial $P_m(S)$ has the expected properties. One checks that
$P_2(S)=S^2+2$. Then $P_{2m}(S)$ is a polynomial in $C=S^2+2$. In fact
\begin{equation}
\sum_{m \geq 0} t^m P_{2m}(S)=\frac{2-tC}{1-tC+t^2}=\sum_{m \geq 0} t^m Q_m(C).
\end{equation}
Comparison of the generating functions shows that 
\begin{equation}
P_m(iS)=i^mQ_m(S) \qquad Q_m(iS)=i^mP_m(S).
\end{equation}
Moreover the definition of the polynomials $P_m$ and $Q_m$ in terms of
$u$ shows that
\begin{equation}
Q_m(C)=P_m(S)^2+2(-1)^{m+1}.
\end{equation}
For analogous reasons, $P_{2m+1}(S)/S$ is also a polynomial in $C$. In fact
\begin{equation}
\sum_{m \geq 0} t^m P_{2m+1}(S)/S = \frac{1+t}{1-tC+t^2}= \sum_{m \geq 0}
t^mR_m(C).
\end{equation}

\section{Relations\label{sect:relations}}

The clue is to consider equation (\ref{eq:scasm}) for $m=l'$.
\begin{equation}\label{eq:scasl'}
\prod_{n=0}^{l'-1}(S-q'^nq^{1/2}k+q'^{-n}q^{-1/2}k^{-1})=\varepsilon(l')
(-\eta)^{l'} f^{l'} e^{l'}.
\end{equation}
The nice feature is that inside the product $q'^n$ runs over all
$l'^{\th}$ roots of unity. The left hand side involves only $S$, $k$
and $k^{-1}$, which commute among themselves. There is a standard way
to simplify the product. We introduce  commuting variables $u$ and $v$
and compute
$$ \begin{array}{rcl}
\prod_{n=0}^{l'-1} (u-u^{-1}-q'^nv+q'^{-n}v^{-1}) & = &
\prod_{n=0}^{l'-1} u^{-1}(u-q'^nv)(u+q'^{-n}v^{-1})  \\ & = &
u^{-l'}\prod_{n=0}^{l'-1} (u-q'^nv) \prod_{n=0}^{l'-1}
(u+q'^{-n}v^{-1}) \\ & = & u^{-l'} (u^{l'}-v^{l'}) (u^{l'}-(-v^
{-1})^{l'}) \\ & = & u^{l'}+(-u^{-1})^{l'}-v^{l'}-(-v^{-1})^{l'}. 
\end{array} $$
In this identity we set $S=u-u^{-1}$ and $v=q^{1/2}k$. On the right hand side
we recognise the Chebychev polynomials from section \ref{sect:functions}. 
Hence (\ref{eq:scasl'}) becomes
\begin{equation} \label{eq:srel}
P_{l'}(S)=q^{l'/2}k^{l'}+(-1)^{l'}q^{-l'/2}k^{-l'}+\varepsilon(l')
(-\eta)^{l'} f^{l'} e^{l'}.
\end{equation}
This gives a relation between $S$ and the $l'^{\th}$ powers of the
other generators. We can now check the properties announced in the
previous section. The right hand side is obtained essentially by
raising to the $l'^{\th}$ power the terms in the equation defining
$S$. The polynomial in $S$ on the left hand side is fixed by its value
in the quotient were $e=f=0$ so that only Cartan generators
 survive.

From this it is clear how to get relations in the centre. We have to
distinguish several cases.

\begin{prop}
  If $l$ is not twice an odd integer, $l'=L$ is even, and (\ref{eq:srel})
  is a relation in the centre :
  \begin{equation} 
    (-1)^{L/2}Q_{L/2}(C)=-k^{L}-k^{-L}+
    \eta^{L} f^{L} e^{L}.
  \end{equation}
\end{prop}

As observed in section
\ref{sect:restrictions}, if $l$ is twice an odd integer, then
$l=L=2l'$ and $l'$ is odd. Then $e^{l'}$ and $f^{l'}$ anti-commute
among themselves and with $k$, $k^{-1}$ and $S$. Hence (\ref{eq:srel})
is a srelation in the scentre. Multiplying by $k^{l'}$ on both
sides gives the relation in the centre written below (\ref{eq:rel2}).
Another relation in the centre (\ref{eq:rel3}) expresses 
$Q_{L/2}(C)=P_{l'}(S)^2+2$ in terms of central elements. 
Using (\ref{eq:rel2})
this gives a formula for $f^{L} e^{L}$ in the standard basis of the
centre (\ref{eq:rel4}). To summarise 

\begin{prop}
   If $l$ is twice an odd integer, the following
  relations in the centre hold:  
  \begin{eqnarray}
    && Sk^{L/2}R_{\frac{L-2}{4}}(C)  =  q^{L/4}k^L-q^{-L/4}+
    (-1)^{\frac{L+2}{4}}\eta^{L/2} f^{L/2} e^{L/2}k^{L/2} \;,
    \label{eq:rel2}
    \\
    && Q_{L/2}(C)  =  
    -k^{L}-k^{-L}+2(-1)^{\frac{L+2}{4}}q^{L/4}\eta^{L/2} f^{L/2}
    e^{L/2}\left(k^{L/2}+k^{-L/2}\right)-\eta^{L} f^{L} e^{L} \;,
    \qquad \label{eq:rel3}
    \\
    && Q_{L/2}(C)  =  k^{L}+k^{-L}+4+2q^{L/4}S\left(k^{L/2}+k^{-L/2}\right)
    R_{\frac{L-2}{4}}(C)-\eta^{L} f^{L} e^{L}.
    \label{eq:rel4}
  \end{eqnarray}
\end{prop}

\section{Representations\label{sect:representations}}

The representation theory of $\cU_q(osp(1|2))$ has already been studied by
several authors \cite{Kobayashi,PSosp,GeSunXue}. It seems to us,
however, that a complete classification of them, including both
periodic and nilpotent ones did not exist in the case $q$ a root
of unity. Our point is also to illustrate the use of the Scasimir in
such a classification, as in \cite{Pinc} for the classical case. 

First we list some families of representations of $\cU_q(osp(1|2))$,
show that they are irreducible and give the possible isomorphisms
among them. After that we show that any irreducible representation
appears in our list.

Let $V$ be a vector space with basis $|0\rangle,\cdots,|L-1\rangle$, 
with the convention that $|L\rangle=|0\rangle$.  We define operators $Q,U,P$
acting on $V$ by
\begin{equation}\label{eq:qup}
  Q|m\rangle=q^{-m}|m\rangle \qquad U|m\rangle=
  (-1)^m|m\rangle \qquad P|m\rangle=|m+1\rangle.
\end{equation}
We endow $V$ with several structures of $\cU_q(osp(1|2))$-modules. 
\begin{itemize}
\item The module $M_+(\lambda,\phi,\sigma)$ depends on three complex
  parameters, the first two are nonzero. As a vector space
  $M_+(\lambda,\phi,\sigma)$ is $V$. We set 
  \begin{equation} \label{eq:dqfe}
    k=\lambda Q \qquad f=\phi P \qquad e=\frac{1}{\eta
      \phi}P^{-1}\left(q^{1/2}\lambda Q-q^{-1/2}\lambda^{-1}
    Q^{-1}-\sigma U\right).  
    \label{eq:M+}
  \end{equation}
  By definition, $\sigma U$ anti-commutes with $f$ and $e$, but commutes
  with $k$, and the definition of $e$ gives it the status of the
  Scasimir operator, so the structure equations of $\cU_q(osp(1|2))$ are
  trivially satisfied.
  
  Conjugation by $Q$, $P$ and $U$ shows that 
  \begin{equation} \label{eq:M+equiv} M_+(\lambda,\phi,\sigma) 
    \cong M_+(\lambda,q^{-1}\phi,\sigma)\cong M_+(q\lambda,\phi,-\sigma)\cong 
    M_+(\lambda,-\phi,\sigma).
  \end{equation}
  
  The value of the central element $e^L$ can be computed in terms of the
  parameters of $M_+(\lambda,\phi,\sigma)$. We write this value as
  $\overline{\epsilon}^L$ for a certain $\overline{\epsilon}$. 
\item The module $M_-(\lambda,\epsilon,\sigma)$ depends on three complex
  parameters, the first two are non\-zero. As a vector space
  $M_-(\lambda,\epsilon,\sigma)$ is $V$. We set 
  \begin{equation} \label{eq:dqef}
    k=\lambda Q \qquad e=\epsilon P^{-1} \qquad f=\frac{1}{\eta
      \epsilon}\left(q^{1/2}\lambda Q-q^{-1/2}\lambda^{-1}
    Q^{-1}-\sigma U \right)P.  
    \label{eq:M-}
  \end{equation}
  Again, for analogous reasons, the structure equations of
  $\cU_q(osp(1|2))$ are 
  trivially satisfied, and there are equivalences. 
\end{itemize}
If we choose the parameters of $M_+(\lambda,\phi,\sigma)$
and $M_-(\lambda,\epsilon,\sigma)$ such that $\overline{\epsilon}=
\epsilon$, then 
$M_+(\lambda,\phi,\sigma)$ $ \cong
M_-(\lambda,\epsilon,\sigma)$.
This is because we can change basis in $M_+(\lambda,\phi,\sigma)$ by
setting $\overline{|m\rangle}=\epsilon^m e^{-m}|m\rangle$. In this new
basis, $e,k,S$ act like on $M_-(\lambda,\epsilon,\sigma)$, and there is
no freedom to define $f$.

We claim that unless $l$ is odd and $\sigma=0$ the above
representations are irreducible. This is because any submodule would
contain a common eigenvector of $k$ and $S$. Assume we are dealing
with a module of type $M_+$ for example. The action of the invertible
operator $f$ will create $L$ non-zero vectors which are distinguished
by the eigenvalues of $k$ and $S$, hence are linearly independent. 

In case $l$ is odd and $\sigma = 0$, the module is irreducible if
considered as graded-module, for the same reasons as above (the
gradation playing the role played by $S$ when $\sigma\neq 0$). 

In contrast with the classical case \cite{Pinc} and the case $q$ is
not a root of unity, an ungraded finite dimensional simple 
module cannot always been endowed
with a gradation, as shown by the following. 
In case $l$ is odd and $\sigma = 0$, the (ungraded) 
module is not irreducible. 
Assume again we are dealing with a module of type $M_+$. The
invertible operator $f^l$ commutes with $P$ and $Q$ (or $k$), hence
also with $e$, but is not scalar. A non trivial eigenspace of $f^l$
is a subrepresentation. We describe explicitly but without details 
the corresponding representations. Consider a vector space $V'$ with
basis $|0\rangle,\cdots,|l-1\rangle$, with the convention that 
$|l\rangle=|0\rangle$.  We define operators $Q,P$ (but not $U$) acting on $V'$ 
by formul\ae\ analogous to (\ref{eq:qup}). Then
\begin{itemize}
\item The module $M_+(\lambda,\phi)$ depends on two non-zero complex
parameters. As a vector space $M_+(\lambda,\phi)$ is $V'$. The action
of $k$, $f$ and $e$ is defined by formul\ae\ (\ref{eq:dqfe}) with
$\sigma=0$.
\item The module $M_-(\lambda,\epsilon)$ depends on two non-zero complex
parameters. As a vector space $M_-(\lambda,\epsilon)$ is $V'$. The action
of $k$, $e$ and $f$ is defined by formul\ae\ (\ref{eq:dqef}) with
$\sigma=0$.
\end{itemize}
The equivalences between those modules follow from the previous
equivalences. Those modules are irreducible.

Now, we describe another type of representations $M_d(\lambda)$ 
where both $e$ and
$f$ are nilpotent. 
Let $V_d$ be a
vector space with basis $|0\rangle,\cdots,|d-1\rangle$, with the
convention that $|-1\rangle=|d\rangle=0$.
We define operators $Q$, $U$ and $P$ as in (\ref{eq:qup}) but note
that $P$ is not invertible anymore ($P|d-1\rangle=0$). We also define
$P'$ by $P'|m\rangle=|m-1\rangle$.
Then we set
\begin{equation} 
k=\lambda Q \qquad f=P \qquad e=\frac{1}{\eta}P'\left(q^{1/2}\lambda
Q-q^{-1/2}\lambda^{-1}Q^{-1}-\sigma U\right).  
\end{equation}
For the time being, $\lambda$ and $\sigma$ are arbitrary complex parameters.
This time, one has to check explicitly the structure equations of 
$\cU_q(osp(1|2))$ because $f=P$ is not invertible, so $P'$ is not an
inverse. The only nontrivial check is that
$ef+fe=\frac{k-k^{-1}}{q-q^{-1}}$ at the boundaries 
$|0\rangle$ and $|d-1\rangle$. One finds
\begin{equation}
q^{1/2}\lambda-q^{-1/2}\lambda^{-1}=\sigma=q^{1/2}q'^{-d}\lambda-q^{-1/2}q'^{d}
\lambda^{-1}.
\end{equation}
These equations give $\sigma$ in terms of $\lambda$ and
\begin{equation}
  (q'^d -1)(\lambda^2 - q'^{d-1}) = 0 \;.
  \label{eq:quantisation}
\end{equation}
Unless $q'^d=1$
(that is unless $d$ is a multiple of $l'$), this is a quantisation
condition for $\lambda$. The presence of
$q'=-q$ (instead of $q$)  in (\ref{eq:quantisation})
explains why the even dimensional modules have no classical limit. 

Now we have to check irreducibility. By construction the basis vectors
of $V_d$ are eigenvectors of $k$ and $S$ and $|d-1\rangle$ is
annihilated by $f$. From section \ref{sect:foundations} we know the
structure of $\cU_q(osp(1|2))$ as an $\cA$-module. It implies that
acting with powers of $e$ on $|d-1\rangle$ generates a submodule. So
the representation can be irreducible only if $e^{d-1}$ does not act
as $0$. The converse is also true because $e$ is nilpotent on $V_d$,
so has a kernel in any submodule. Hence $V_d$ is irreducible if and
only if $e^{d-1}|d-1\rangle\neq 0$, which is equivalent to
$f^{d-1}e^{d-1}|d-1\rangle\neq 0$. This last condition is easy to
check using (\ref{eq:scasm}) for $m=d-1$. A simple computation gives
the irreducibility criterion :
\begin{equation}
  \prod_{n=1}^{d-1}
  (q'^n-1)(q^{1/2-d}\lambda+q^{d-1/2}\lambda^{-1}q'^{-n})\neq 0.
  \label{eq:irred}
\end{equation}
Equivalently, $V_d$ is irreducible if either $d<l'$ or $d=l'$ and
$q'\lambda^2$ is not a nontrivial power of~$q'$. 

It is now easy to show 
\begin{theo}
  The finite dimensional irreducible representations of $\cU_q(osp(1|2))$
  are:
  \begin{itemize}
  \item The $f$-periodic modules $M_+(\lambda,\phi,\sigma)$ 
    (\ref{eq:M+}) of
    dimension $L$ (take instead $M_+(\lambda,\phi)$ 
    of dimension $l$ 
    if $l$ is odd and $\sigma=0$)
  \item The $e$-periodic modules $M_-(\lambda,\epsilon,\sigma)$
    (\ref{eq:M-}) of
    dimension $L$ (take instead $M_-(\lambda,\epsilon)$ of dimension $l$ 
    if $l$ is odd and $\sigma=0$)
  \item The nilpotent modules $M_d(\lambda)$, of dimension 
    $1\le d\le l'$, with the conditions
    (\ref{eq:quantisation}) and (\ref{eq:irred}). 
  \end{itemize}
  Modules of type $M_+$ and $M_-$ are equivalent if, and only if, they
  share the action of the central elements. 
\end{theo}
{\em Proof:} 
What remains to prove is that every
irreducible finite dimensional representation of $\cU_q(osp(1|2))$ 
appears in the above list

The proof goes as follows.  
Let $V$ be an irreducible finite dimensional representation space of
$\cU_q(osp(1|2))$. The operators $k$ and $S$ can be simultaneously
diagonalised on $V$. The reason is that they commute, so they have a
common eigenvector. As a consequence of the commutations relations of
the generators of $\cU_q(osp(1|2))$ with $k$ and $S$, monomials in 
$e,f,k,k^{-1}$ applied to this vector are clearly either $0$ or again 
eigenstates $k$ and $S$. Those monomials span a subrepresentation
which has to be $V$ itself by irreducibility. So we endow $V$ with a
basis consisting of eigenvectors of $k$ and $S$.

Now we distinguish three cases. First, suppose $f$ is invertible on
$V$. Then knowing the action of $k$, $S$ and $f$ fixes the action of
$e$ uniquely by $\eta e=f^{-1}(q^{1/2}k-q^{-1/2}k^{-1}-S)$. As $f^L$
is central, starting from a common eigenvector of $k$ and $S$ and
acting with $f$ one builds a subrepresentation of dimension at most
$L$, which has to be $V$ itself. Unless $l$ is twice an odd integer
and $S=0$, $V$ has to be $L$-dimensional because $1,f,\cdots,f^{L-1}$
are distinguished by the eigenvalue of $k$ and $S$. If $l$ is twice an
odd integer and $S=0$ we can diagonalise $f^{l'}$ and $k$. So we
always end-up with an irreducible representation of type $M_+$.
The second case is when $e$ is invertible. The same line of arguments
leads to an irreducible representation of type $M_-$. The third case
is when neither $f$ nor $e$ is invertible. We start from an
eigenvector of $k$ and $S$ annihilated by $e$ and call it the highest
weight vector. The structure of $\cU_q(osp(1|2))$ as an $\cA$-module
implies that powers of $f$ acting on the highest weight build a
subrepresentation which has to be $V$ itself. The largest
non-vanishing power, say $d-1$, of $f$ acting on the highest weight gives a
lowest weight. For the same reasons, powers of $e$ acting on the
lowest weight build a subrepresentation which has to be $V$
itself. The formula relating $k,f,e$ to $S$ shows that $V$ has to
be equivalent to a representation on $V_d$ listed above. We have
already given a criterion for irreducibility for those.

\qed


\newpage

\begin{thebibliography}{99}

\bibitem {RosHC} M. Rosso, 
  {\sl Analogues de la forme de Killing et du th{\'e}or{\`e}me
    d'Harish--Chandra pour les groupes quantiques,} 
  Ann. Scient. {\'E}c. Norm. Sup., $4^e$ s{\'e}rie, t.23, (1990), 445.
  
\bibitem {Lusztigmod} G. Lusztig, 
  {\sl Modular representations and quantum groups,}
  Contemp. Math. {\bf 82} (1989) 59.
  
\bibitem {DK}  C. De Concini and V.G. Kac, 
  {\sl Representations of quantum groups at roots of 1,}
  Progress in Math. {\bf 92} (1990) 471 (Birkh{\"a}user).
  
\bibitem {KerDipl} T. Kerler,
  {\sl Darstellungen der Quantengruppen und Anwendungen,}
  Diplomarbeit, ETH--Zurich, August 1989.

\bibitem {ABrelations} D. Arnaudon and M. Bauer,
  {\sl Polynomial Relations in the Centre of 
    ${\cal U}_q(sl(N))$, } 
  hep-th/9310030, 
  Lett. Math. Phys. {\bf 30} (1994) 251.
  
\bibitem {AACuqsln} B. Abdesselam, D. Arnaudon and A. Chakrabarti,
  {\sl Representations of ${\cal U}_q(sl(N))$ at Roots of Unity,}
  q-alg/9504006,   
  J. Phys. A: Math. Gen. {\bf 28} (1995) 5495.

\bibitem{KacSuperRep} V.G. Kac, 
  {\sl Representations of classical Lie
    superalgebras}, in Lecture Notes in Mathematics {\bf 676},
  Springer-Verlag, Berlin, Heidelberg, New York, 1978.
  
\bibitem {AAB} B. Abdesselam, D. Arnaudon and M. Bauer, {\sl Centre
    and Representations of ${\cal U}_{q}(sl(2|1))$ at Roots of Unity,}
  q-alg/9605015, 
  ENSLAPP-A-583/96.

\bibitem {KulRes} P.P. Kulish and N. Yu Reshetikhin, 
  {\sl Universal $R$-matrix of the quantum superalgebra $osp(1|2)$,} 
  Lett. Math. Phys. {\bf 18} (1989).

\bibitem {Lesniewski} A. Le\'sniewski,
  {\sl A remark on the Casimir elements of Lie superalgebras and
    quantized Lie superalgebras,}
  J. Math. Phys. {\bf 36} (3) (1995) 1457.

\bibitem {PaisRitt} A. Pais and V. Rittenberg,
  {\sl Semisimple graded Lie algebras,}
  Journ. Math. Phys. {\bf 16} (1975) 2062.

\bibitem {Pinc} G. Pinczon, 
  {\sl The enveloping algebra of the Lie  superalgebra $osp(1|2)$,}
  Journ. Algebra {\bf 132} (1990) 219.

\bibitem {ABFscasimirs} D. Arnaudon, M. Bauer and L. Frappat, {\sl On
    Casimir's ghost,} q-alg/9605021, ENSLAPP-A-587/96.

\bibitem {Duflo} M. Duflo, 
  {\sl Sur la classification des id{\'e}aux primitifs dans l'alg{\`e}bre
    enveloppante d'une alg{\`e}bre de Lie semi-simple,}
  Ann. of Math. (2) {\bf 105} (1977) 107.

\bibitem {MussonPrim} I.M. Musson,
  {\sl A classification of primitive ideals in the enveloping algebra
    of a classical simple Lie superalgebra,}
  Adv. in Math. {\bf 91} (1992) 252.

\bibitem {Zhangd} R. B. Zhang, 
  {\sl Finite-dimensional representations of $U_q(osp(1|2n))$ and its
    connection with quantum $so(2n+1)$,}
  Lett. Math. Phys. {\bf 25} (1992) 317.

\bibitem {Kassel} C. Kassel,
  {\sl Quantum Groups,} 
  GTM {\bf 155}, Springer Verlag 1995. 
 
\bibitem {Kobayashi} K.-I. Kobayashi, {\sl Representation theory of
    $osp(1|2)_q$,} Z. Phys. {\bf C 59} (1993) 155.

\bibitem {PSosp} T. D. Palev and N. I. Stoilova, {\sl Unitarizable
    representations of the deformed parabose superalgebra
    $U_q(osp(1/2))$ at roots of 1,} 
  q-alg/9507026, 
  J. Phys. {\bf A28} (1995) 7275.

\bibitem {GeSunXue} Ge Mo-Lin, Sun Chang-Pu and Xue Kang,
  {\sl New $R$-matrices for the Yang--Baxter equation associated with
    the representations of the quantum superalgebra $U_q(osp(1,2))$
    with $q$ a root of unity,}
  Phys. Lett. {\bf A 163} (1992) 176.

\end{thebibliography}
\end{document}